\documentclass[12pt]{iopart}

\usepackage{amssymb}
\usepackage{graphicx,epsfig}
\def\bpn{\bigskip\par\noindent}

\def\.{{\;}}

\def\req#1{{\rm(\ref{eq:#1})}}

\def\ve{{\varepsilon}}

\def\sgn{\mbox{sgn}}

\newcommand\R{\mathbb{R}}

\newtheorem{theorem}{Theorem}[section]

\newtheorem{lemma}[theorem]{Lemma}

\newtheorem{proposition}[theorem]{Proposition}

\newtheorem{remark}[theorem]{Remark}

\begin{document}

\title[Constraint regularization]
{On the Relation between Constraint Regularization, Level Sets, and
Shape Optimization}

\author{A.~Leit\~ao\dag\ and O.~Scherzer\ddag
\footnote[3]{Correspondence should be sent to otmar.scherzer@uibk.ac.at} }

\address{\dag\ Department of Mathematics,
           Federal University of St.\,Catarina, PO Box 476, 88010-970
           Florianopolis, Brazil, e-mail: {\tt aleitao@mtm.ufsc.br} }

\address{\ddag\ Department of Computer Science,
           University Innsbruck, Technikerstra{\ss}e 25, A-6020
           Innsbruck, Austria, e-mail: {\tt Otmar.Scherzer@uibk.ac.at} }

\begin{abstract}
We consider regularization methods based on the coupling of
Tikhonov regularization and projection strategies.
From the resulting constraint regularization method we
obtain level set methods in a straight forward way.
Moreover, we show that this approach links the areas of
asymptotic regularization to inverse problems theory,
scale-space theory to computer vision, level set methods,
and shape optimization.
\end{abstract}

\section{Introduction}

The major goal of this paper is to highlight the relation between the
following areas:
\begin{enumerate}
\item Regularization for inverse and ill-posed problems, in particular
      \begin{enumerate}
      \item Tikhonov regularization for constraint operator equations
      \item Asymptotic regularization
      \end{enumerate}
\item Scale-space theory in \emph{computer vision}
\item Shape optimization
\end{enumerate}

The general context is to solve the constraint ill-posed operator
equation:
\begin{equation}
\label{eq:op}
F(u)=y\,,
\end{equation}
where $u$ is in the admissible class
\begin{equation*}
U:=\{u : u = P(\phi) \mbox{ and } \phi \in {\cal{D}}(P)\}\;.
\end{equation*}

The constraint equation can be
formulated as an unconstrained equation
\begin{equation}
\label{eq:op-u}
F(P(\phi)) = y\;.
\end{equation}
Assuming that the operator equation is ill-posed it has to be regularized
for a stable solution.

Classical results on convergence and stability of regularization
(see e.g. \cite{Mor84,Mor93,EngHanNeu96}) such as
\begin{enumerate}
\item existence of a regularized solution
\item stability of the regularized approximations
\item approximation properties of the regularized solutions
\end{enumerate}
are applicable if $P$ is
\begin{enumerate}
\item bounded and linear or
\item nonlinear, continuous, and weakly closed.
\end{enumerate}
In order to link constraint regularization methods, shape optimization,
level sets, and inverse scale-space, we require discontinuous operators $P$,
and thus the classical framework of regularization theory is not applicable
yet.

Tikhonov regularization for solving the unconstrained equation \req{op}
consists in approximation the solution of \req{op} by the minimizer
$u_\alpha$ of the functional
\begin{equation*}
\|F(u)-y\|^2 + \alpha \|u-u_*\|^2\;.
\end{equation*}
If $F$ is differentiable, then
\begin{equation}
\label{eq:opta}
F'(u_\alpha)^*(F(u_\alpha)-y) + \alpha (u_\alpha-u_*) = 0\,,
\end{equation}
where $F'(u_\alpha)^*$ denotes the adjoint $\cdot^*$ of the
derivative of $f$ at $u_\alpha$.
\req{opta} is the optimality condition for a minimizer of the
Tikhonov functional.
Using the formal setting $\Delta t := 1/\alpha$,
$u(\Delta t):=u_\alpha$, and $u(0):=u_*$ we find
\begin{equation*}
F'(u(\Delta t))^*(F(u(\Delta t))-y) +
\frac{u(\Delta t)-u(0)}{\Delta t} = 0\;.
\end{equation*}
Thus $u_\alpha=u(\Delta t)$ can be considered as the solution of
one implicit time step with step-length $\Delta t = \frac{1}{\alpha}$
for solving
\begin{equation}
\label{eq:asmp}
\frac{\partial u}{\partial t} = - F'(u)^*(F(u)-y)
\end{equation}
we end up with the \emph{inverse scale-space method} (see e.g.
\cite{GroSch00,SchGro01}).
We note that the inverse scale-space method corresponds to the
\emph{asymptotic regularization method} as introduced by Tautenhahn
\cite{Tau94,Tau95}.

The terminology "inverse scale-space" is motivated from
\emph{scale-space} theory in \emph{computer vision}: images contain structures
at a variety of scales. Any feature can optimally be recognized at
a particular scale. If the optimal scale is not available a-priori,
it is desirable to have an image representation at
multiple scales.

A \emph{scale-space} \index{scale-space} is an image representation at
a continuum of scales, embedding the image $u$ into a family
\begin{equation*}
\{ T_t (u) : t \geq 0 \}
\end{equation*}
of gradually simplified versions satisfying:
\begin{enumerate}
\item \emph{Recursivit\"at:}
      \begin{equation*}
      T_0(u) = u\;.
      \end{equation*}
\item \emph{Kausalit\"at:}
      \begin{equation*}
      T_{t+s}(u) = T_t(T_s(u)) \mbox{ for all } s,t \geq 0\;.
      \end{equation*}
\item \emph{Regularit\"at:}
      \begin{equation*}
      \lim_{t \to 0+} T_{t}(u) = u\;.
      \end{equation*}
\end{enumerate}
For more background on the topic of scale-space theory we refer to
\cite{Lin94,NieJohOlsWei99,Wei98a,Ker01}.

The ill--posedness of inverse problems prohibits such a representation in
scales of images and the concept has to be replaced by inverse scale-space
theory, which includes approximative causality together with:
\begin{enumerate}
\item \emph{Inverse Rekursivit\"at:}
      \begin{equation*}
      T_\infty(y) = u^\dagger\;.
      \end{equation*}
\item \emph{Inverse Regularit\"at:}
      \begin{equation*}
      \lim_{t \to \infty-} T_{t}(y) = u^\dagger\;.
      \end{equation*}
\end{enumerate}
Here $y$ is the input data and $u^\dagger$ is a solution
of \req{op}.
As shown in \cite{SchGro01} \req{asmp} is an inverse scale
space method.

In this work we show that the inverse scale-space method for the
constrained inverse problem \req{op-u} with appropriate $P$ is a
\emph{level set method}.
Level set methods have been developed by Osher \& Sethian \cite{OshSet88}
(see also Sethian \cite{Set99}).
Recently, level set methods have been successfully applied for the solution
of inverse problems (see, e.g.,
\cite{San96,LitLesSan98,LitLesSan98,DorMilRap00,ItoKunLi02,RamLamLesZol01,RamLamLes01,Bur01}).

Moreover, we show that the shape derivative in form optimization
and the level set derivative correspond. For simplicity of
presentation we concentrate on highlighting this link by
considering a particular example from \cite{HetRun96}.

\section{Derivation of the Level Set Method}

In this section we consider the constraint optimization problem
of solving \req{op} on the set of piecewise constant functions which
attain two values, which we fix for the sake of simplicity of
presentation to $0$ and $1$. Typical examples include parameter
identification problems where the value $1$ denotes an inclusion.

Let $\Omega \subseteq \R^n$ be bounded with boundary $\partial \Omega$
Lipschitz.
Set
\begin{equation*}
{\cal{P}} := \{u : u = \chi_{\tilde{\Omega}}:
     \tilde{\Omega} \subseteq \Omega \} \cap L^2(\Omega)\,,
\end{equation*}
then the unconstrained inverse problem consists in solving \req{op-u} with
\begin{equation*}
\begin{array}{rl}
P : H^1(\Omega) &\to {\cal{P}}\;.\\
     \phi &\mapsto
           \frac{1}{2} + \frac{1}{2} \sgn (\phi) =:
           \frac{1}{2} + \frac{1}{2} \left\{ \begin{array}{rl}
                                 1 \mbox{ for } \phi \geq 0\\
                                 -1 \mbox{ for } \phi < 0
                                 \end{array} \right.
\end{array}
\end{equation*}
Moreover, let for the sake of simplicity of presentation,
\begin{equation*}
F: L^2(\Omega) \to L^2(\Omega)
\end{equation*}
be Fr\'echet-differentiable. It is as well possible to consider the
operator $F$ in various Hilbert space settings such as for instance
$F:H^1(\Omega) \to L^2(\partial \Omega)$. Since it does not make any
methodological differences we concentrate on an operator on $L^2(\Omega)$.
Also the space $H^1(\Omega)$ is chosen more or less arbitrarily; we have
selected these spaces in such a way that the typical distance functions for smooth
domains are contained in $H^1(\Omega)$.
\bpn
Tikhonov regularization for this problem consists in minimizing the
functional
\begin{equation}
\label{eq:unperturbed}
\int_\Omega (F(P(\phi)) - y)^2 + \alpha
\int_\Omega \left((\phi - \phi_*)^2 + |\nabla (\phi-\phi_*)|^2\right) \;.
\end{equation}
Since the functional does not attain a minimum, we consider the 
``minimizer'' $\phi_\alpha$ as
\begin{equation*}
\phi_\alpha = \lim_{\ve \to 0+} \phi_{\ve,\alpha}\,,
\end{equation*}
where $\phi_{\ve,\alpha}$ minimizes the functional
\begin{equation}
\label{eq:ve}
\int_\Omega (F(P_\ve(\phi)) - y)^2 + \alpha
\int_\Omega \left((\phi - \phi_*)^2 + |\nabla (\phi-\phi_*)|^2\right) \;.
\end{equation}
We use
\begin{equation*}
P_\ve(t) := \left\{
\begin{array}{rcl}
0 & \mbox{ for } & t < -\ve\,,\\
1 + \frac{t}{\ve} & \mbox{ for } & t \in \left[-\ve,0\right]\,,\\
1 &\mbox{ for } & t > 0\,,
\end{array} \right.
\end{equation*}
for approximating $P$ as $\ve \to 0^+$.
In this case we have
\begin{equation*}
P'(t) = \lim_{\ve \to 0+} P'_\ve(t) = \delta(t)\;.
\end{equation*}
Here and in the following $\delta(t)$ denotes the one-dimensional
$\delta$-distribution.
Moreover, we denote
\begin{equation*}
u_\alpha := \lim_{\ve \to 0+} P_\ve(\phi_{\alpha,\ve})\;.
\end{equation*}
Note that we do \emph{not} require that $u_\alpha = P(\phi_\alpha)$.
The proposed methodology to define generalized solutions
$u_\alpha = \lim_{\ve \to 0+} P(\phi_{\ve,\alpha})$ is a standard way 
in \emph{phase transitions}.

In the following we derive an optimality condition for a minimizer of
\req{unperturbed}, which is considered the limit $\ve \to 0+$ of the
minimizers of the functionals \req{ve}. For this purpose it is
convenient to recall some basic results from \emph{Morse theory} of surfaces.
The particular results are collected from \cite{FomKun97}. We emphasize
that in this paper we only apply the Morse theory to compact, smooth subset
of $\R^2$, which of course can be considered as surfaces.
\begin{proposition}
\label{pr:1}
Let $\phi$ be a smooth function on a compact smooth surface $M$, and
$\phi^{-1}[a,b] \subseteq M$ contain no critical point of $\phi$.
Then,
\begin{enumerate}
\item the level sets $\phi^{-1}(b)$ and $\phi^{-1}(a)$ are diffeomorphic
      (in particular they consist of the same number of smooth circles
       diffeomorphic to a standard circle)
      \cite[Proposition 6.2.1.]{FomKun97}.
      In particular the Hausdorff measure of
      $\phi^{-1}(t), t \in [a,b]$ changes continuously.
\item Moreover, for any $\rho \in [a,b]$, $\phi^{-1}(\rho)$ is a smooth
      compact
      $1$-manifold \cite[p107]{FomKun97}.
      In particular $\phi^{-1}(\rho)$ can be parameterized by
      finitely many disjoint curves.
\end{enumerate}

\end{proposition}

\begin{lemma}
\label{le:important}
Let $\phi$ be a smooth function, having no critical point in a compact
neighborhood $M$ of the level set $\phi^{-1}(0)$.
Then,
\begin{equation*}
\lim_{\ve \to 0+} P_\ve'(\phi) =
\frac{1}{|\nabla \phi|} \delta(\phi)\;.
\end{equation*}
We recall that $\delta(\phi)$ is the one-dimensional
$\delta$-distribution centered at the level line in normal direction.
\end{lemma}
{\bf{Proof:}}
In dimension $1$ this is a well-known result, especially in
physics (see \cite{SchwQM1,SchwQM2}). 
We sketch the proof adopted
to level set functions in dimension $2$; for higher dimension the
generalization is obvious.

From Proposition \ref{pr:1} we know that the level set
$\phi^{-1}(0)$ is a smooth compact $1$-manifold, which can be
parameterized by a curve $s(\tau), \tau \in [0,2\pi)$
\footnote[1]{For the sake of simplicity of presentation we assume
that the level set is parameterized by just one curve. The general
case of finitely many disjoint curves is analogous.}, i.e.,
\begin{equation*}
\phi^{-1}(0) := \left\{ s(\tau)=(s_1(\tau),s_2(\tau)) :
 \tau \in [0,2\pi) \right\}\;.
\end{equation*}

Here $n$ is the normal vector to the level set, which can be
characterized as
\begin{equation*}
n(\tau) = -\frac{\nabla \phi}{|\nabla \phi|}(s(\tau))\;.
\end{equation*}
We choose the negative sign in the definition of the normal vector
based on the following considerations: if $\phi$ is a monotonically
increasing function in normal direction to the level set pointing
into the domain bounded by the level set, then
$n(\tau)$, as defined above, points outside this domain.

The basic idea of the proof is to find a relation between a parameter
$\ve$ and a parametric function $\psi:[0,2\pi) \to \R$ such that the sets
\begin{equation*}
\Omega_\psi := \left\{ s(\tau) + \rho n(\tau) : \tau \in [0,2\pi),
                       \rho \in [0,\psi(\tau))\right\}
\end{equation*}
and $\phi^{-1}(-\ve,0]$ ``asymptotically'' correspond.
\bpn
By making a Taylor series expansion we find
\begin{equation*}
~\hspace{-2.5cm}
\begin{array}{rcl}
\displaystyle
\phi(\Omega_\psi) &=&
\displaystyle
\phi \left(\left\{ s(\tau) - \rho  \frac{\nabla \phi}{|\nabla \phi|}(s(\tau))
                  : \tau \in [0,2\pi),
                  \rho \in [0,\psi(\tau))\right\} \right)\\
&=& \displaystyle
    \left\{ \phi \left(
                 s(\tau) - \rho  \frac{\nabla \phi}{|\nabla \phi|}(s(\tau))
                 \right) :
                 \tau \in [0,2\pi), \rho \in [0,\psi(\tau)) \right\}\\
&=& \displaystyle
    \left\{ \phi (s(\tau))  -
            \rho \frac{\nabla \phi}{|\nabla \phi|}(s(\tau))
                 \nabla \phi (s(\tau)) + O(\rho^2):
                 \tau \in [0,2\pi), \rho \in [0,\psi(\tau)) \right\}\\
&=& \displaystyle
    \left\{ - \rho |\nabla \phi|(s(\tau)) + O(\rho^2):
                 \tau \in [0,2\pi), \rho \in [0,\psi(\tau)) \right\}\;.
\end{array}
\end{equation*}
If we choose
\begin{equation*}
\psi(\tau) := \psi_\ve (\tau) = \frac{\ve}{|\nabla \phi (s(\tau))|}\,,
\end{equation*}
and set
\begin{equation*}
C_{\min} := \inf \{ |\nabla \phi|(s(\tau)): \tau \in [0,2\pi)\}\,,
\end{equation*}
then there exists a constant $C$ such that
\begin{equation*}
\hspace*{-2cm}
\displaystyle
\Omega_- :=
\left[ - \ve + \ve^2 \frac{C}{C_{\mbox{min}}^2},-\ve^2
                     \frac{C}{C_{\mbox{min}}^2} \right]
\subseteq
\phi(\Omega_\psi) \subseteq
\left[ - \ve - \ve^2 \frac{C}{C_{\mbox{min}}^2},\ve^2 \frac{C}{C_{\mbox{min}}^2}
\right] =: \Omega_+\;.
\end{equation*}
Set $\tau = \frac{C}{C_{\mbox{min}}^2}$. Then, for
$v \in C(\overline{\Omega})$, it follows from the
\emph{coarea formula} \cite{EvaGar92} that
\begin{equation*}
~\hspace{-1.5cm}
\begin{array}{rl}
\displaystyle
~ &\left| \int_{\phi^{-1}(-\ve,0)} v - \int_{\Omega_\psi} v \right|\\
\leq &
\displaystyle
\frac{\mbox{ max } |v|}{C_{\mbox{min}}}
\left\{ \int_{\phi^{-1}(-\ve-\tau\ve^2,-\ve+\tau\ve^2)} |\nabla \phi| +
        \int_{\phi^{-1}(-\tau\ve^2,\tau\ve^2)} |\nabla \phi| \right\}\\
\leq&
\displaystyle
\frac{\mbox{ max } |v|}{C_{\mbox{min}}}
\left\{ \int_{-\ve-\tau\ve^2}^{-\ve+\tau\ve^2}
        {\cal{H}}^{1}(\phi^{-1}(\rho)) d\rho +
        \int_{-\tau\ve^2}^{\tau\ve^2}
        {\cal{H}}^{1}(\phi^{-1}(\rho)) d\rho \right\}
\end{array}
\end{equation*}
where ${\cal{H}}^{1}(\phi^{-1}(\rho))$ is the $1$-dimensional Hausdorff
measure of the set $\phi^{-1}(\rho)$. According to Proposition \ref{pr:1}
${\cal{H}}^{1}(\phi^{-1}(\rho))$ is uniformly bounded.
This implies that
\begin{equation*}
\displaystyle
\left| \int_{\phi^{-1}(-\ve,0)} v - \int_{\Omega_\psi} v \right|=O(\ve^2)\,,
\end{equation*}
and consequently
\begin{equation*}
~\hspace{-1.5cm}
\begin{array}{rcl}
\displaystyle
\lim_{\ve \to 0+} \int_\Omega P_\ve'(\phi) v &=&
\displaystyle
\lim_{\ve \to 0+} \frac{1}{\ve} \int_{\phi^{-1}(-\ve,0)} v\\
&=&
\displaystyle\lim_{\ve \to 0+}
\frac{1}{\ve} \int_{\Omega_{\psi_\ve}} v\\
&=&
\displaystyle \lim_{\ve \to 0+} \int_{\Omega_{\psi_\ve}}
\frac{1}{\psi_\ve} \frac{1}{|\nabla \phi|}  v\;.
\end{array}
\end{equation*}
This shows that
\begin{equation*}
~\hspace*{-2.5cm}
\begin{array}{rcl}
\displaystyle
\lim_{\ve \to 0+} \int_\Omega P_\ve'(\phi) v
&=&
\displaystyle \lim_{\psi \to 0+}
\int_0^{2\pi} \frac{1}{|\nabla \phi|(s(\tau))} \frac{1}{\psi(\tau)} \cdot \\
& & \displaystyle \qquad  \quad \cdot
              \int_0^{\psi} v \left| \det
              \left[ \begin{array}{rl}
                     s_1'(\tau)+\rho n_1'(\tau) & n_1(\tau)\\
                     s_2'(\tau)+\rho n_2'(\tau) & n_2(\tau)
                     \end{array}
              \right] \right| d \rho d\tau\\
&=& \displaystyle
    \int_\Omega \delta(\phi) \frac{v}{|\nabla \phi|}\;.
\end{array}
\end{equation*}
~\hfill $\Box$
\bpn
Lemma \ref{le:important} is central to derive the optimality condition
for a minimizer of \req{unperturbed}.

From the definition of a minimizer of \req{ve} it follows that
for all $h \in H^1(\Omega)$
\begin{equation}
\label{eq:opt-weak}
\begin{array}{rl}
~ & \displaystyle
    \int_\Omega (F(u_{\ve,\alpha})-y) F'(u_{\ve,\alpha})
P_\ve'(\phi_{\ve,\alpha}) h \\
& \displaystyle
 + \alpha \int_\Omega \left( (\phi_{\ve,\alpha}-\phi_*)h +
                          \nabla(\phi_{\ve,\alpha}-\phi_*) \nabla h
                   \right) = 0\;.
\end{array}
\end{equation}
We denote by $F'(u)^*,P_\ve'(\phi)^*$ the $L^2$-adjoints of $F'(u)$,
$P_\ve'(\phi)$, respectively, i.e., for all $v,w \in L^2(\Omega)$
\begin{equation*}
\int_\Omega w (F'(u) v) = \int_\Omega (F'(u)^*w) v \mbox{ and }
\int_\Omega w (P_\ve'(\phi) v) = \int_\Omega (P_\ve'(\phi)^*w) v\;.
\end{equation*}
Since $P_\ve'(\phi)$ is self-adjoint, i.e., $P_\ve'(\phi)^* = P_\ve'(\phi)\,,$
it follows that
\begin{equation}
\label{eq:opt}
\begin{array}{rl}
\displaystyle
P_\ve'(\phi_{\ve,\alpha}) F'(u_{\ve,\alpha})^* (F(u_{\ve,\alpha})-y) +
\alpha (I-\Delta) (\phi_{\ve,\alpha}-\phi_*) &= 0
\mbox{ on } \Omega\,,\\
\displaystyle
\frac{\partial (\phi_{\ve,\alpha} -\phi_*)}{\partial n} &=0 \mbox{ on } \partial \Omega\;.
\end{array}
\end{equation}
Thus $u_\alpha=\lim_{\ve \to 0+}u_{\ve,\alpha}$ and
$\phi_\alpha=\lim_{\ve \to 0+}\phi_{\ve,\alpha}$ satisfies
\begin{equation}
\label{eq:reg}
\delta (\phi_\alpha)
\frac{F'(u_\alpha)^* (F(u_\alpha)-y)}{|\nabla \phi_\alpha|} +
\alpha (I-\Delta) (\phi_\alpha-\phi_*) = 0\;.
\end{equation}
For the sake of simplicity of presentation we assume that the
operator $F$ is of such quality that $F'(u)^* (F(u)-y)$ is continuous
on $\Omega$. Note that in general this may not be the case since
$F'(u)^* (F(u)-y) \in H^1(\Omega).$

Therefore, it follows from \req{reg} that
\begin{equation*}
(I - \Delta)^{-1}
\left(\delta(\phi_\alpha)
\frac{F'(u_\alpha)^* (F(u_\alpha)-y)}{|\nabla \phi_\alpha|}
\right) +
\alpha (\phi_\alpha-\phi_*) = 0\;.
\end{equation*}
Set $\alpha=\frac{1}{\Delta t}$ and set $\phi_\alpha = \phi(t)$,
$\phi_* = \phi(0)$ and accordingly $u(t)=P(\phi(t))$.
Then, by taking the formal limit $\Delta t \to 0+$ we get
the asymptotic regularization method
\begin{equation}
\label{eq:pde-weak}
\frac{\partial \phi}{\partial t} =
- (I - \Delta)^{-1}
\left( \delta(\phi(t))
       \frac{F'(u(t))^*(F(u(t)) -y)}{|\nabla \phi(t)|}
  \right)\;.
\end{equation}
The right hand side $v$ of \req{pde-weak} solves the equation
\begin{equation}
\label{eq:potential}
\begin{array}{rcl}
(I - \Delta) v &=& - \delta(\phi(t))
\frac{F'(u(t))^*(F(u(t)) -y)}{|\nabla \phi(t)|} \,,\\
\frac{\partial v}{\partial n} &=& 0\;.
\end{array}
\end{equation}
Using potential theory (see e.g. \cite{Kre99,Eng97})
a solution $v_1$ of the homogeneous problem
\begin{equation*}
\Delta v_1 (t) =
\delta(\phi(t))
\frac{F'(u(t))^*(F(u(t)) -y)}{|\nabla \phi (t)|}\,,
\end{equation*}
is given by the \emph{single layer potential}
\begin{equation*}
v_1(x) = -
\int_{\phi(t)^{-1}(0)}
      \frac{F'(u(t))^*(F(u(t)) -y)(z) \gamma(x,z)}
           {|\nabla \phi(t)(z)|} \,dz\,,
\end{equation*}
where
\begin{equation}
\label{eq:12}
\gamma(x,y) = \left\{
\begin{array}{rl}
\frac{1}{2\pi} \ln \left( \frac{1}{|x-y|} \right) & \mbox{ in } \R^2\,,\\
\frac{1}{4\pi} \frac{1}{|x-y|} & \mbox{ in } \R^3
\end{array} \right.
\end{equation}
is the \emph{single layer potential}.

Then, $v=v_1+v_2$ solves \req{potential} where $v_2$ solves
\begin{equation*}
\begin{array}{rl}
v_2-\Delta v_2 &=-v_1 \mbox{ on } \Omega\\
\frac{\partial v_2}{\partial n} &= -\frac{\partial v_1}{\partial n}
\mbox{ on } \partial \Omega\;.
\end{array}
\end{equation*}
\req{pde-weak} is a \emph{level set method} describing the evolution 
of the level set function $\phi$. The zero level set of $\phi$, i.e., the set
$\{\phi=0\}$, describes the boundary of the inclusions to be
recovered.

\begin{remark}
An adequate approximation of $P$ is central in our considerations.
The family of functions
\begin{equation*}
Q_\ve(t) := \left\{
\begin{array}{rcl}
0 & \mbox{ for } & t < -\ve\,,\\
\frac{t+\ve}{2\ve} & \mbox{ for } & t \in \left[-\ve,\ve\right]\,,\\
1 &\mbox{ for } & t > \ve\,,
\end{array} \right.
\end{equation*}
approximates the $\delta$-distribution too.
Since the point-wise limit of $Q_\ve$ is
\begin{equation*}
P(t) := \left\{ \begin{array}{rl}
                     0 & \mbox{ for } t < 0\,,\\
                     \frac{1}{2} & \mbox{ for } t=0\,,\\
                     1 & \mbox{ for } t > 0\,,
             \end{array} \right.
\end{equation*}
which is not in ${\cal{P}}$ if the $n$-dimensional Lebesgue measure
of $\phi^{-1}(0)$ is greater than zero. This would not be appropriate
for our problem setting.
\end{remark}

In this section we have elaborated on the interaction between constraint
regularization methods and level set methods. We have shown that our level
set method can be considered as an inverse scale-space method,
respectively asymptotic regularization method.
In contrast to standard results on asymptotic regularization methods and
inverse scale-space methods (see \cite{Tau94,Tau95,GroSch00}), here
the situation is more involved, since the regularizer of the
underlying regularization functional \req{unperturbed} is considered
as approximation of the minimizers of the functional \req{ve}, i.e.,
it is a $\Gamma$-limit (see e.g. \cite{AmbDan99}).

One of the most significant advantages of level set methods is
that the topology of the zero--level set may change over time. So
far, this situation has not covered by our derivation of level set
methods, where we essentially relied on Proposition \ref{pr:1} and
Lemma \ref{le:important}. In case a topology change occurs the
Morse index of the level set function $\phi$ changes and
Proposition \ref{pr:1} and consequently Lemma \ref{le:important}
are not applicable. Moreover, in this case the single layer
potential representations \req{12} are no longer valid 
(see e.g. \cite{ColKre83,Kre99}), since the topology changes results 
in domain with cusps. The effect of topology changes on the level set 
methods are status of ongoing research. 
In this article we are interested in revealing
interactions between constraint regularization techniques, level
set methods, and shape optimization. To show the interaction part
we rely on some explicit calculations of the shape derivative in
\cite{HetRun96} where inclusions are considered smooth without
cusps. Thus in order to compare level set evolution and shape
derivative, we find it desirable to limit our considerations and 
neglect topology changes.

\subsection{Relation to other level set methods}

\req{pde-weak} is a Hamilton-Jacobi type equation of the form
\begin{equation}
\label{eq:level2}
\frac{\partial \phi}{\partial t} + V \nabla \phi = 0
\end{equation}
with velocity
\begin{equation*}
V =
\frac{(I-\Delta)^{-1}
\left( \delta(\phi) \frac{F'(u(t))^*(F(u(t)) -y)}{|\nabla \phi (t)|}
       \right)}{|\nabla \phi (t)|}
\frac{\nabla \phi}{|\nabla \phi (t)|}\;.
\end{equation*}

The numerical solution of \req{pde-weak} is similar to the implementation
of well-established level set methods, like e.g. considered
by Santosa \cite{San96}, who suggested a velocity
\begin{equation*}
V = - F'(u(t))^*(F(u(t)) -y) \frac{\nabla \phi (t)}{|\nabla \phi (t)|}\;.
\end{equation*}
The differential equation
\begin{equation}
\label{eq:level}
\frac{\partial \phi}{\partial t} =
F'(u(t))^*(F(u(t)) -y) |\nabla \phi(t)|
\end{equation}
is solved explicit in time, which results in
\begin{equation*}
\frac{\phi(t+\Delta t) - \phi(t)}{\Delta t} =
F'(u(t))^*(F(u(t)) -y) |\nabla \phi(t)|\;.
\end{equation*}
After several numerical time-steps the iterates are \emph{updated}.
In our level-set approach such an update is inherent, since in each
step the data is normalized by the operator $(I-\Delta)^{-1}$.

\subsection{Relation to Shape Optimization}

In this subsection we show that the term
\begin{equation*}
\delta(\phi) \frac{F'(u)^* (F(u)-y)}{|\nabla \phi|}
\end{equation*}
is the steepest descent direction of the functional
$\|F(u)-y\|^2$ with respect to the \emph{shape} of the
level set $\phi^{-1}(0)$.

It is much more illustrative to show this relation exemplary. To this
end we consider the \emph{inverse potential problem} of recovery of
a object $D \subseteq \R^2$ in
\begin{equation*}
\Delta v = \chi (D) \mbox{ in }\Omega \mbox{ with } v=0 \mbox{ on }
\partial \Omega\;.
\end{equation*}
In this context
\begin{equation*}
\begin{array}{rcl}
F : L^2(\Omega) & \to & L^2(\Omega)\;.\\
    f & \mapsto & \Delta^{-1} f \mbox{ with homogeneous Dirichlet data}
\end{array}
\end{equation*}
The numerical recovery of shape of the inclusion $D$ from Neumann boundary
measurements was considered in \cite{HetRun96}.
For the sake of simplicity of presentation, here we are interested in
the shape derivative of $F$, while Hettlich and Rundell
considered the operator $T \circ F$, where $T$ is the Neumann trace operator.
Since $T$ is linear the shape derivative of $T \circ F$ is completely
determined by the shape derivative of $F$, and thus we do not impose
any restriction on the consideration by considering the simpler problem.

The operator $F$ is linear and thus the Gateaux-derivative of $F$ at $u$
in direction $h$ satisfies $F'(u) h = F(h)$.
Thus the \emph{level set derivative} is given by
\begin{equation}
\label{eq:no}
v:= F'(u) P'(\phi) h = F(P'(\phi) h) =
\Delta^{-1} \left( \delta(\phi) \frac{h}{|\nabla \phi|} \right)\;.
\end{equation}
Let $v_1$ be the single layer potential according to $h$ on $\phi^{-1}(0)$,
i.e.,
\begin{equation*}
v_1(x) = - \int_{\phi^{-1}(0)} \frac{1}{2\pi}
\ln \frac{1}{|x-y|} \frac{h}{|\nabla \phi|}(y)\,dy\;.
\end{equation*}
This function satisfies
\begin{equation*}
\Delta v_1 = \delta(\phi) \frac{h}{|\nabla \phi|} \mbox{ on }
\Omega\;.
\end{equation*}
Let $v_2$ be the solution of
\begin{equation*}
\Delta v_2 = 0 \mbox{ on } \Omega \mbox{ and } v_1=-v_2 \mbox{ on }
\partial \Omega\;.
\end{equation*}
Then $v = v_1 + v_2$ solves
\begin{equation*}
\Delta v = \delta(\phi) \frac{h}{|\nabla \phi|} \mbox{ on } \Omega
\mbox{ and } v=0 \mbox{ on }
\partial \Omega\;.
\end{equation*}
Moreover, the single layer potential satisfies on the zero level set
\begin{equation*}
\begin{array}{rcl}
\left(\frac{\partial v_1}{\partial n} \right)_+ -
\left(\frac{\partial v_1}{\partial n} \right)_- &=&
\frac{h}{|\nabla \phi|}\,,\\
(v_1)_+&=&(v_1)_-\;.
\end{array}
\end{equation*}
Here $(\cdot)_+$, $(\cdot)_-$ denote the limits from outside, inside
of the domain bounded by the zero level curves, respectively.

We recall that $h$ is considered a perturbation of the level set
\emph{function}. A change in the level set function implies a change
in the zero level set, which eventually turns out to be the shape derivative.

To make this concrete, let $s_{th}$ the parameterizations of
$(\phi+th)^{-1}(0)$, i.e., \\
$(\phi+th)(s_{th})=0$. We make a Taylor Ansatz with respect to the
parametrization
\begin{equation}
\label{eq:expansion}
s_{th} = s + t \tilde{h} + O(t^2)\,,
\end{equation}
and a series expansion for $\phi$ and $h$, which gives
\begin{equation*}
0 = (\phi+th)(s_{th})
  = t \nabla \phi \tilde{h} + t h(s) + O(t^2)\;.
\end{equation*}
This shows that on the zero level set we have
\begin{equation*}
\frac{h}{|\nabla \phi|} = -
\frac{\nabla \phi}{|\nabla \phi|} \cdot \tilde{h} =
n \cdot \tilde{h}\;.
\end{equation*}

Thus $v$ satisfies the differential equation
\begin{equation}
\label{eq:shape}
\begin{array}{c}
\left\{ \begin{array}{rcl}
        \displaystyle
        \Delta v &=&
        \displaystyle 0 \mbox{  on } \Omega \backslash \phi^{-1}(0)\,,\\
        \displaystyle v &=&
        \displaystyle 0 \mbox{  on } \partial \Omega\;;
    \end{array}\right.\\
~ \\
\displaystyle
\left(\frac{\partial v}{\partial n} \right)_+ -
\left(\frac{\partial v}{\partial n} \right)_- =
\tilde{h} \cdot n \mbox{ on } \phi^{-1}(0)\,,\\
~\\
\displaystyle
(v)_+ = (v)_-  \mbox{ on } \phi^{-1}(0)\;.
\end{array}
\end{equation}
This is the shape derivative $F'(D)(\tilde{h})$ of $F$ at
$D=\{x : P(\phi) > 0\}$ in direction $\tilde{h}$ as calculated by
Hettlich and Rundell \cite{HetRun96}.

Our calculations show the level set derivative $v:=F'(u)P'(\phi)h$
can be computed from the shape derivative. Now, we point out
that the converse is evenly true. This is nontrivial since the
arguments $\tilde{h}$ appearing in the shape derivative are
multidimensional functions, while the argument $h$
in the level set derivative is one-dimensional.

Let $\tilde{h}$ be expressed in terms of the local
coordinate system $n$ and $\tau$, where $n$, $\tau$ are the normal,
respectively tangential vectors on the zero level set, i.e.,
\begin{equation*}
\tilde{h} = h n + h_\tau \tau\;.
\end{equation*}
The shape derivative is independent of the tangential component, which
in particular implies that the shape derivative gradient descent
deforms the shapes in normal direction to the level curve.
Thus, from \req{no} we find that
\begin{equation}
\label{eq:relation}
F'(D)(\tilde{h}) = F'(D)(h n) = F'(P \phi) h\;.
\end{equation}
That is we have shown:
\begin{theorem}
\label{th:relation}
By \req{relation} the level set derivative $F'(u)P'(\phi)h = F(P'(\phi))h$
is uniquely determined from the shape derivative and vice versa.
\end{theorem}
From Theorem \ref{th:relation} we see that the level set derivative
moves the zero level set in direction of the shape derivative.

\end{document}